# Deformation from symmetry for Schrödinger equations of higher order on unbounded domains[*]

Addolorata SALVATORE
Dipartimento Interuniversitario di Matematica
Via Orabona 4, I–70123 Bari, Italy
e-mail: salvator@dm.uniba.it

Marco SQUASSINA
Dipartimento di Matematica e Fisica
Via Musei 41, I–25121 Brescia, Italy
e-mail: m.squassina@dmf.unicatt.it

March 20, 2018

## Abstract

By means of a perturbation method recently introduced in [7], we discuss the existence of infinitely many solutions for a class of perturbed symmetric higher order Schrödinger equations with non–homogeneous boundary data on unbounded domains.

## 1 Introduction and main results

Let $\Omega$ be an unbounded domain in $\mathbb{R}^N$ with $\partial\Omega$ smooth, $N \geq 2K$ and $K \geq 1$, $\varphi$ and $V$ two functions in suitable $L^s$ spaces. The main goal of this note is to study the existence of multiple solutions for the following polyharmonic Schrödinger equation

$$\begin{cases} (-\Delta)^K u + V(x)u = g(x,u) + \varphi & \text{in } \Omega \\ \left(\frac{\partial}{\partial \nu}\right)^j u = \phi_j & \text{on } \partial\Omega \\ D^j u(x) \to 0 & \text{as } |x| \to \infty \\ j = 0, \dots, K-1 \end{cases} \quad (1)$$

where $\nu$ is the unit outward normal to $\partial\Omega$ and $\phi_j$ belongs to $H^{K-j-\frac{1}{2}}(\partial\Omega)$ for $j = 0, \dots, K-1$.

So far, many papers have been written on the existence and multiplicity of solutions for *second order* elliptic problems with Dirichlet boundary data, especially by means of variational methods. In particular, problem (1) with $K = 1$, $\phi_0 = 0$, $V = 0$ and $\Omega$ *bounded* has been studied by many researchers in the last decades. If $\varphi = 0$ and $g(x, \cdot)$ is odd, the problem is symmetric and multiplicity results can be proven in a standard fashion for *any* subcritical $g$ (see [20] and references therein). On the contrary, if $\varphi \neq 0$ the symmetry is lost and a natural question is whether the multiplicity is preserved under perturbation of $g$. Partial answers have been given in [2, 3, 4, 14, 19, 21] where existence of infinitely many solutions was obtained via techniques of critical point theory, provided that suitable restrictions on the growth rate of $g$ are assumed. Roughly speaking, if $g(x, u) \simeq |u|^{p-2}u$ the exponent $p$ is required to be greater than 2 but not too large, that is $2 < p < 2 + \frac{2}{N-2}$.

[*]This work was partially supported by Ministero dell'Università e della Ricerca Scientifica e Tecnologica (40% – 2001) and by Gruppo Nazionale per l'Analisi Funzionale e le sue Applicazioni.

The success in looking for solutions of these problems made quite interesting to study the case where, in general, the boundary datum $\phi_0 \in H^{1/2}(\partial\Omega)$ may differs from zero. This introduces a higher order loss of symmetry, since the functional associated with the problems contains two terms which fail to be even. In this situation some multiplicity achievements have been proved in [8, 9, 10, 11, 12] under suitable stronger restrictions on the growth of $g$ and on the regularity of $\Omega$; in particular, if $g(x, u) \simeq |u|^{p-2}u$, infinitely many solutions have been found if $2 < p < 2 + \frac{2}{N-1}$. Finally, for unbounded domains and $V \neq 0$, some results have been recently established by one of the authors in [17, 18] under more involved assumptions on the growth of $g$.

Now, a natural question is whether these results for the second order case extend to the *higher order*. If $\Omega$ is bounded, some results have been recently obtained in [13]; in this case, the "critical" exponent for the problem becomes $2 + \frac{2K}{N}$ if $\phi_j \neq 0$ for some $j$ while it is $2 + \frac{2K}{N-2K}$ (the natural extension of $2 + \frac{2}{N-2}$ to the case $K > 1$) if $\phi_j = 0$ for any $j$. In the case where $\Omega$ is *unbounded* and $K > 1$ no multiplicity result for (1) is, to our knowledge, known. We will show (see Theorems 1.1 and 1.2) that results similar to the second order case hold; in particular, if the Lebesgue measure of $\Omega$ is finite, we find again the results of [13]. In Corollary 1.3, we give a multiplicity result when $\Omega = \mathbb{R}^N$. We recall that, in a variational setting, if $\Omega$ is unbounded, these problems also present a lack of compactness due to the failure of compactness of the Sobolev embedding $H^K(\Omega) \hookrightarrow L^2(\Omega)$. In order to overcome this problem we assume, besides suitable assumptions on $G$, that the function $V$ has a "good" behaviour at infinity so that the Schrödinger operator $(-\Delta)^K + V$ on $L^2(\Omega)$ admits a discrete spectrum (this may fail, in general, if $V(x) \to 0$ as $|x| \to \infty$) and the Palais–Smale condition can be recovered. If $V = 0$ our condition reduces to a "shrinking" assumption on $\Omega$ at infinity which again allow us to regain compactness and also the discreteness of the spectrum of $(-\Delta)^K$. Notice that this last property may fail on general domains of $\mathbb{R}^N$. If for instance $\Omega$ is connected with smooth boundary such that

$$\{x \in \mathbb{R}^N : x_1 = \cdots = x_{N-1} = 0, \ x_N > 0\} \subset \Omega \subset \{x \in \mathbb{R}^N : x_N > 0\}$$

and $\partial x_N / \partial \nu \leq 0$, then $-\Delta$ admits a purely continuous spectrum (see [15]).

We will actually obtain multiplicity results for a class of higher order operators more general than $(-\Delta)^K$. To achieve this we apply a method recently developed by Bolle, Ghoussoub and Tehrani for dealing with problems with perturbed symmetry (see [7, 8] for the abstract framework and [12] for some recent generalizations).

To state the main results, we assume that the following conditions hold:

($g_1$) $g : \Omega \times \mathbb{R} \to \mathbb{R}$ is a continuous map and $g(x, \cdot)$ is odd;

($g_2$) there exists $\bar{s} \neq 0$ such that $\inf_{\Omega} G(x, \bar{s}) > 0$;

($g_3$) there exists $\mu > 2$ such that for every $x \in \Omega$, $s \in \mathbb{R}$, $s \neq 0$

$$0 < \mu G(x,s) \leq s g(x,s);$$

($g_4$) there exist $c > 0$ and $2 < p < \mu+1$, $p < K_*$ if $N > 2K$, such that for every $(x,s) \in \Omega \times \mathbb{R}$

$$|g(x,s)| \leq c|s|^{p-1}$$

where $K_* = \frac{2N}{N-2K}$ is the critical exponent for the Sobolev embedding of $H_0^K(\Omega)$;

(A) $A$ is a formally selfadjoint elliptic differential operator of order $2K$ with constant coefficients and there exists $\gamma > 0$ such that

$$\forall u \in C_c^\infty(\Omega): \quad \int_\Omega Au \cdot u \geq \gamma \begin{cases} \int_\Omega |\Delta^m u|^2 & \text{if } K = 2m \\ \int_\Omega |\nabla \Delta^m u|^2 & \text{if } K = 2m+1; \end{cases} \quad (2)$$

(V) $V \in L^2_{loc}(\Omega)$ is such that $\underset{\Omega}{\mathrm{infess}}\, V > 0$ and

$$\lim_{|x| \to \infty} \int_{S(x) \cap \Omega} \frac{1}{V(\xi)} d\xi = 0$$

where $S(x)$ is the unit ball of $\mathbb{R}^N$ centered at $x$.

Assumption (V) has been used by Benci and Fortunato in [5] for proving some compact embedding theorems for weighted Sobolev spaces. It is easy to see that (V) holds in particular if $V$ is a continuous positive function on $\mathbb{R}^N$ which goes to infinity as $|x| \to \infty$. As we will see, this assumption implies that the spectrum of $A + V(x)$ with Dirichlet boundary conditions in $L^2(\Omega)$ is discrete (see Proposition 2.1); from now on we will denote by $(\lambda_n)$ the divergent sequence of its eigenvalues (counted with their multiplicity). We denote by $L^2(\Omega, V)$ a suitable weighted $L^2$ space (see Section 2 for further details) and by $\mu'$ the conjugate exponent of $\mu$. Moreover, we write $|B|$ for the Lebesgue measure of the set $B$ in $\mathbb{R}^N$.

The following are the main results in the case $V \neq 0$.

**Theorem 1.1.** *Assume that conditions $(g_1) - (g_4)$, (A) and (V) hold. Let $\varphi \in L^{\mu'}(\Omega)$ and*

$$\phi_j \in H^{K-j-\frac{1}{2}}(\partial \Omega) \qquad j = 0, \ldots, K-1.$$

*Moreover, suppose that for $n \in \mathbb{N}$ sufficiently large we have*

$$\lambda_n > n^{\kappa(p,K,N,\mu)}, \qquad \kappa(p,K,N,\mu) = \frac{2K\mu(p-2)}{(\mu-p+1)(2Kp-(p-2)N)}. \quad (3)$$

*Then, the boundary value problem*

$$\begin{cases} Au + V(x)u = g(x,u) + \varphi & \text{in } \Omega \\ \left(\frac{\partial}{\partial \nu}\right)^j u = \phi_j & \text{on } \partial\Omega \\ D^j u(x) \to 0 & \text{as } |x| \to \infty \\ j = 0, \ldots, K-1 \end{cases} \quad (S_{\varphi,\phi})$$

admits an unbounded sequence of solutions $(u_n) \subset H^K(\Omega) \cap L^2(\Omega, V)$. Moreover, the same conclusion holds provided that in place of (3) we have

$$\frac{\mu}{\mu - p + 1} < \frac{2Kp}{N(p-2)}$$

under the additional assumption that $|\Omega| < \infty$.

If $\phi_j = 0$ for every $j = 0, \ldots, K-1$, then the previous result can be improved.

**Theorem 1.2.** *Assume that conditions $(g_1) - (g_4)$, $(A)$ and $(V)$ hold. Moreover, let $\varphi \in L^{\mu'}(\Omega)$ and suppose that for $n \in \mathbb{N}$ sufficienly large we have*

$$\lambda_n > n^{\overline{\kappa}(p,K,N,\mu)}, \qquad \overline{\kappa}(p, K, N, \mu) = \frac{2K\mu(p-2)}{(\mu - 1)(2Kp - (p-2)N)}. \tag{4}$$

*Then, the boundary value problem*

$$\begin{cases} Au + V(x)u = g(x, u) + \varphi & \text{in } \Omega \\ \left(\frac{\partial}{\partial \nu}\right)^j u = 0 & \text{on } \partial \Omega \\ D^j u(x) \to 0 & \text{as } |x| \to \infty \\ j = 0, \ldots, K-1 \end{cases} \tag{$S_{\varphi,0}$}$$

*admits an unbounded sequence of solutions $(u_n) \subset H_0^K(\Omega) \cap L^2(\Omega, V)$. Moreover, the same conclusion holds provided that in place of (4) we have*

$$\frac{\mu}{\mu - 1} < \frac{2Kp}{N(p-2)}$$

*under the additional assumption that $|\Omega| < \infty$.*

In general, conditions (3) and (4) are verified if $V$ has a fast growth at infinity and $p$ is greater than 2 but not too large. Notice that in the case $|\Omega| < \infty$ the potential does *not* affect the multiplicity range anymore. We now give an application of Theorem 1.2 when $\Omega = \mathbb{R}^N$ by considering a more particular class of potentials $V$ such that the growth of the eigenvalues $(\lambda_n)$ of $A + V(x)$ can be explicitely estimated.

**Corollary 1.3.** *Assume that conditions $(g_1)-(g_4)$, $(A)$ and $(V)$ hold with $\Omega = \mathbb{R}^N$. Moreover, assume that there exists $\alpha \geq 1$ such that*

$$\limsup_{\lambda \to +\infty} \frac{|\{x \in \mathbb{R}^N : V(x) < \lambda\}|}{\lambda^{\frac{N}{\alpha}}} < +\infty. \tag{5}$$

*Then, for every $\varphi \in L^{\mu'}(\mathbb{R}^N)$ the problem*

$$\begin{cases} Au + V(x)u = g(x, u) + \varphi & \text{in } \mathbb{R}^N \\ D^j u(x) \to 0 & \text{as } |x| \to \infty \\ j = 0, \ldots, K-1 \end{cases} \tag{6}$$

*admits un unbounded sequence of solutions* $(u_n) \subset H^K(\mathbb{R}^N) \cap L^2(\mathbb{R}^N, V)$ *provided that*

$$\frac{\alpha}{2K + \alpha} > \frac{\mu(p-2)N}{(\mu-1)(2Kp - (p-2)N)}.$$

*In particular, if* $\mu = p$, *then* (6) *admits infinitely many solutions for any* $p \in ]2, \bar{p}[$ *where* $\bar{p}$ *is the largest root of the quadratic equation*

$$2(KN + \alpha(N - K))p^2 - (\alpha(5N - 2K) + 4KN)p + 2\alpha N = 0.$$

*Proof.* Denote by $\mathcal{N}(\lambda, A + V(x), \mathbb{R}^N)$ the number of the eigenvalues of $A + V(x)$ in $L^2(\mathbb{R}^N)$ which are less or equal than $\lambda$. As proved in [16, Theorem 3], there exists a constant $B_{N,K} > 0$ such that

$$\mathcal{N}(\lambda, A + V(x), \mathbb{R}^N) \leq B_{N,K} \int_{\mathbb{R}^N} \left((\lambda - V(x))_+\right)^{N/2K}$$

for every $\lambda$. Clearly, by virtue of the positivity of $V$ and (5), for $\lambda$ sufficiently large we have

$$\int_{\mathbb{R}^N} \left((\lambda - V(x))_+\right)^{N/2K} = \int_{\{V(x) < \lambda\}} (\lambda - V(x))^{N/2K}$$
$$\leq \lambda^{N/2K} \left|\{V(x) < \lambda/2\}\right| \leq M \lambda^{\frac{N(2K+\alpha)}{2K\alpha}}$$

for some positive constant $M$ depending on $N$ and $\alpha$. Therefore, by choosing $\lambda = n$ in the previous inequality, for $n \in \mathbb{N}$ sufficiently large we have

$$\lambda_n \geq C_{N,K,\alpha} \, n^{\frac{2K\alpha}{N(2K+\alpha)}}$$

being $C_{N,K,\alpha}$ a suitable positive constant depending on $N, K$ and $\alpha$. The assertion now follows immediately by applying Theorem 1.2. □

In order to give an idea of the amplitude of the range $]2, \bar{p}[$, in the following table we list the values of $\bar{p}$ for $K = 1, \ldots, 10$ corresponding to the dimensions $N = 2K + 1$.

| $\alpha \geq 1$ | $K = 1$ | $K = 2$ | $K = 3$ | $\cdots$ | $K = 10$ |
|---|---|---|---|---|---|
| $\alpha = 1$ | $\bar{p} = 2.2310$ | $\bar{p} = 2.1688$ | $\bar{p} = 2.1284$ | $\cdots$ | $\bar{p} = 2.0463$ |
| $\alpha = 2$ | $\bar{p} = 2.3494$ | $\bar{p} = 2.2895$ | $\bar{p} = 2.2319$ | $\cdots$ | $\bar{p} = 2.0901$ |
| $\alpha = 3$ | $\bar{p} = 2.4201$ | $\bar{p} = 2.3786$ | $\bar{p} = 2.3161$ | $\cdots$ | $\bar{p} = 2.1314$ |
| $\alpha = 4$ | $\bar{p} = 2.4668$ | $\bar{p} = 2.4466$ | $\bar{p} = 2.3854$ | $\cdots$ | $\bar{p} = 2.1704$ |
| $\alpha = 5$ | $\bar{p} = 2.5000$ | $\bar{p} = 2.5000$ | $\bar{p} = 2.4432$ | $\cdots$ | $\bar{p} = 2.2072$ |
| $\vdots$ | $\vdots$ | $\vdots$ | $\vdots$ | $\ddots$ | $\vdots$ |
| $\alpha = \infty$ | $\bar{p} \simeq 2.6929$ | $\bar{p} \simeq 2.9314$ | $\bar{p} \simeq 3.0514$ | $\cdots$ | $\bar{p} \simeq 3.2816$ |

Table 1: Values of $\bar{p}$ varying $\alpha$ when $K = 1, \ldots, 10$ and $N = 2K + 1$.

Observe that condition (5) holds for example if $V$ is a positive function verifying $V(x) = |x|^\alpha$ for some $\alpha \geq 1$ and for all $x \in \mathbb{R}^N$, $|x|$ large. If, in particular, $V$ grows exponentially fast and $\mu = p$, then Corollary 1.3 yields infinitely many solutions for any $p \in ]2, p_\infty[$, being $p_\infty$ the largest root of the equation $2(N-K)p^2 - (5N - 2K)p + 2N = 0$ (see the last row of Table 1). If $N = 2K + 1$, we get $p_\infty \to \sqrt{2} + 2$ as $K \to \infty$.

Notice that the condition $(V)$ holds if and only if for every $b > 0$

$$\lim_{|x| \to \infty} |S(x) \cap \Omega_b| = 0 \quad \text{where} \quad \Omega_b = \{x \in \Omega : V(x) \leq b\}. \tag{7}$$

In particular, in the case where the potential function $V$ is identically equal to zero, (7) reduces to the following assumption, which corresponds to a "shrinking" condition at infinity:

(D) $\Omega$ is an unbounded domain such that

$$\lim_{|x| \to \infty} |S(x) \cap \Omega| = 0.$$

Condition $(D)$ is a necessary and sufficient condition for the embedding of a space "like" $H^K(\Omega)$ in $L^2(\Omega)$ to be compact; moreover it implies that the spectrum of $A$ with Dirichlet boundary data consists of a sequence $(\mu_n)$ of eigenvalues (with finite multiplicity) having $+\infty$ as the only accumulation point (see Proposition 2.2).

The following corollaries complement the results of [13] dealing with bounded domains.

**Corollary 1.4.** *Assume that $\Omega$ is a domain satisfying $(D)$ and that conditions $(g_1) - (g_4)$ and $(A)$ hold. Let $\varphi \in L^{\mu'}(\Omega)$ and*

$$\phi_j \in H^{K-j-\frac{1}{2}}(\partial \Omega) \qquad j = 0, \ldots, K-1.$$

*Moreover, suppose that for $n \in \mathbb{N}$ sufficienly large we have*

$$\mu_n > n^{\kappa(p,K,N,\mu)}, \qquad \kappa(p,K,N,\mu) = \frac{2K\mu(p-2)}{(\mu-p+1)(2Kp-(p-2)N)}. \tag{8}$$

*Then, the boundary value problem*

$$\begin{cases} Au = g(x,u) + \varphi & \text{in } \Omega \\ \left(\frac{\partial}{\partial \nu}\right)^j u = \phi_j & \text{on } \partial \Omega \\ D^j u(x) \to 0 & \text{as } |x| \to \infty \\ j = 0, \ldots, K-1 \end{cases} \qquad (Q_{\varphi,\phi})$$

*admits an unbounded sequence of solutions $(u_n) \subset H^K(\Omega)$. Moreover, the same conclusion holds provided that in place of (8) we have*

$$\frac{\mu}{\mu - p + 1} < \frac{2Kp}{N(p-2)}$$

*under the additional assumption that $|\Omega| < \infty$.*

If $\phi_j = 0$ for every $j = 0, \ldots, K-1$, a stronger result holds.

**Corollary 1.5.** *Assume that $\Omega$ is a domain satisfying $(D)$. Let $\varphi \in L^{\mu'}(\Omega)$ and assume conditions $(g_1) - (g_4)$ and $(A)$. Moreover, suppose that for $n \in \mathbb{N}$ sufficiently large we have*

$$\mu_n > n^{\overline{\kappa}(p,K,N,\mu)}, \qquad \overline{\kappa}(p,K,N,\mu) = \frac{2K\mu(p-2)}{(\mu-1)(2Kp-(p-2)N)}. \tag{9}$$

*Then, the boundary value problem*

$$\begin{cases} Au = g(x,u) + \varphi & \text{in } \Omega \\ \left(\frac{\partial}{\partial \nu}\right)^j u = 0 & \text{on } \partial\Omega \\ D^j u(x) \to 0 & \text{as } |x| \to \infty \\ j = 0, \ldots, K-1 \end{cases} \qquad (Q_{\varphi,0})$$

*admits an unbounded sequence of solutions* $(u_n) \subset H_0^K(\Omega)$. *Moreover, the same conclusion holds provided that in place of* (9) *we have*

$$\frac{\mu}{\mu - 1} < \frac{2Kp}{N(p-2)}$$

*under the additional assumption that* $|\Omega| < \infty$.

In the case where $\mu = p$ and $\Omega$ has finite measure, we compare in the following table some of the existence ranges for nonzero and zero boundary data when $N = 2K + 1$.

| $K \geq 1$ | $\exists j : \phi_j \neq 0$ | $\forall j : \phi_j = 0$ |
|---|---|---|
| $K = 1$ | $p < 2.6666$ | $p < 4$ |
| $K = 2$ | $p < 2.8000$ | $p < 6$ |
| $K = 3$ | $p < 2.8571$ | $p < 8$ |
| $\vdots$ | $\vdots$ | $\vdots$ |
| $K = 10$ | $p < 2.9523$ | $p < 22$ |

Table 2: Comparison between zero and nonzero boundary data.

Notice, by comparing Table 1 with the second column of Table 2, how the situation $|\Omega| = \infty$ has a *bad* influence, with respect to the case $|\Omega| < \infty$, on the existence ranges. Finally we stress that under suitable *additional* regularity assumptions on $\Omega$, $\phi_j$ and $\varphi$, by virtue of a special Pohožaev type identity, in the case $K = 1$ Theorem 1.1 and Corollary 1.4, can be improved (see [18], in particular Lemma 4.6). Unfortunately, it seems that a similar identity cannot be easily obtained when $K > 1$. On the other hand, if $\phi_j = 0$ for every $j$, the "critical" exponent for our problem seems to be $2 + \frac{2K}{N-2K}$, then the results contained in Theorem 1.2 and Corollary 1.5 coincide with those already stated for $K = 1$ (see [18, Corollary 1.6]) if $|\Omega| = +\infty$ and [3, 8, 21] if $|\Omega| < +\infty$).

## 2 The variational framework

Let $K \geq 1$ and $B \subset \mathbb{R}^N$ smooth. We endow the spaces $L^s(B)$ and $H^K(B)$ with the norms

$$\|u\|_{L^s(B)} = \left(\int_B |u|^s\right)^{1/s}, \quad \|u\|_{H^K(B)} = \left\{\int_B |u|^2 + \sum_{|\mu|=K} \int_B |D^\mu u|^2\right\}^{1/2}.$$

We recall that, by [1, Corollary 4.16], the norm $\|\cdot\|_{H^K(B)}$ is equivalent to the standard norm of $H^K(B)$. Moreover, let $H_0^K(B)$ be the completion of $C_c^\infty(B)$ with respect to $\|\cdot\|_{H^K(B)}$.

Now, we endow $H_0^K(B)$ with another norm equivalent to $\|\cdot\|_{H^K(B)}$. We say that a function $u$ on $B$ is in $\widetilde{H}^K(B)$ if it is the restriction to $B$ of a function in $H^K(\mathbb{R}^N)$. We set

$$\|u\|_{\widetilde{H}^K(B)} = \inf\left\{\|v\|_{H^K(\mathbb{R}^N)}: \ v \in H^K(\mathbb{R}^N),\ v = u \text{ on } B\right\}.$$

It is possible to prove that $\widetilde{H}^K(B)$ is a Banach space, $H_0^K(B)$ is continuously embedded in $\widetilde{H}^K(B)$ and that the norms $\|\cdot\|_{\widetilde{H}^K(B)}$ and $\|\cdot\|_{H^K(B)}$ are equivalent in $H_0^K(B)$ (see [6]). For the sake of simplicity, if $B = \Omega$ we will write $\|u\|_s$, $\|u\|_{K,2}$ and $\|\widetilde{u}\|_{K,2}$ in place of $\|u\|_{L^s(\Omega)}$, $\|u\|_{H^K(\Omega)}$ and $\|u\|_{\widetilde{H}^K(\Omega)}$. From now on, we assume that the function $V$ satisfies condition $(V)$. Then, we can consider the weighted $L^2$ space

$$L^2(\Omega, V) = \left\{u \in L^2(\Omega): \ \int_\Omega V(x)u^2 < +\infty\right\}$$

equipped with the inner product $\int_\Omega V(x)uv$ and the Sobolev space

$$H^K(\Omega, V) = H^K(\Omega) \cap L^2(\Omega, V)$$

endowed with the weighted inner product

$$(u,v)_V = \int_\Omega V(x)uv + \sum_{|\mu|=K} \int_\Omega D^\mu u D^\mu v. \tag{10}$$

We denote by $\|\cdot\|_V$ the norm induced by (10). On the other hand, we can consider the Banach space

$$\widetilde{H}^K(\Omega, V) = \widetilde{H}^K(\Omega) \cap L^2(\Omega, V)$$

endowed with the corresponding norm $\|\cdot\|_{\widetilde{V}} = \|\cdot\|_{K,2} + \int_\Omega V(x)u^2$ and the subspace $H_0^K(\Omega, V) = H_0^K(\Omega) \cap L^2(\Omega, V)$, where $\|\cdot\|_V$ and $\|\cdot\|_{\widetilde{V}}$ are equivalent.

In order to overcome the lack of compactness of the problem, the following Propositions are needed.

**Proposition 2.1.** *Assume that $\Omega$ is an unbounded domain in $\mathbb{R}^N$ and let $V$ a function satisfying assumption $(V)$. Then the embedding $\widetilde{H}^K(\Omega, V) \hookrightarrow L^2(\Omega)$ is compact. It follows that the embedding $H_0^K(\Omega, V) \hookrightarrow L^s(\Omega)$ is compact for every $s \in [2, K_*[$ and the spectrum of the selfadjoint realization of $A + V(x)$ with homogeneous Dirichlet boundary conditions in $L^2(\Omega)$ is discrete.*

*Proof.* If $\Omega = \mathbb{R}^N$, in [5, Theorem 3.1] it has been proved that the space $H^K(\mathbb{R}^N, V)$ is compactly embedded in $L^2(\mathbb{R}^N)$. Small modifications in their proof allow to extend this result to a general unbounded domain $\Omega$ for the space $\widetilde{H}^K(\Omega, V)$. For the sake of completeness we sketch the proof. Let $(u_n)$ be a sequence in $\widetilde{H}^K(\Omega, V)$ such that $u_n \rightharpoonup 0$ in $\widetilde{H}^K(\Omega, V)$. Then, there exists $M \in \mathbb{R}^+$ such that $\int_\Omega V(x)u_n^2 \leq M$ and $u_n|_A \rightharpoonup 0$ in $\widetilde{H}^K(A)$ for every $A \in \sigma_0$, where

$$\sigma_0 = \{A \in \sigma_m: \ A \text{ is open}\},$$
$$\sigma_m = \{A \subset \Omega: \ A \text{ is measurable and } |A \cap S(x)| \to 0 \text{ for } |x| \to \infty\}.$$

By virtue of Theorem 2.8 of [6] we have

$$\forall A \in \sigma_0: \quad \widetilde{H}^K(A) \text{ is compactly embedded in } L^2(A)$$

and then

$$\forall A \in \sigma_0: \quad u_n|_A \to 0 \text{ in } L^2(A). \tag{11}$$

Our aim is to prove that $\|u_n\|_2 \to 0$. Taking any $\varepsilon > 0$, we have

$$\|u_n\|_2^2 = \int_\Omega \frac{1}{V(x)} V(x) u_n^2 \leq \varepsilon M + \int_{\Omega_{1/\varepsilon}} u_n^2,$$

where, by (7), it is $\Omega_{1/\varepsilon} \in \sigma_m$. A slightly modified version of Lemma 3.2 in [5] implies that there exists $A_{1/\varepsilon}$ in $\sigma_0$ such that $\Omega_{1/\varepsilon} \subset A_{1/\varepsilon}$. Hence

$$\forall n \in \mathbb{N}: \quad \|u_n\|_2^2 \leq \varepsilon M + \int_{A_{1/\varepsilon}} u_n^2$$

so that, by (11), $\|u_n\|_2 \to 0$. Therefore, $\widetilde{H}^K(\Omega, V)$ is compactly embedded in $L^2(\Omega)$. Moreover, by the Sobolev embedding we have $H^K(\Omega) \hookrightarrow L^s(\Omega)$ for any $s \in [2, K_*]$. Now, by the Gagliardo–Nirenberg interpolation inequality, for any $u$ which belongs to $L^2(\Omega) \cap L^{K_*}(\Omega)$ it results $u \in L^s(\Omega)$ and $\|u\|_s \leq \|u\|_2^{1-\ell} \|u\|_{K_*}^\ell$ with $\frac{1}{s} = \frac{1-\ell}{2} + \frac{\ell}{K_*}$ and $0 \leq \ell \leq 1$. Hence, the embedding of $H_0^K(\Omega, V)$ in $L^s(\Omega)$ is compact for any $s \in [2, K_*[$. Finally, since the operator $A + V(x)$ with homogeneous Dirichlet boundary conditions is essentially selfadjoint on $C_0^\infty(\Omega)$, the discreteness of the spectrum follows arguing as in [5, Theorem 4.1]. $\square$

**Proposition 2.2.** *Assume that $\Omega$ is an unbounded domain in $\mathbb{R}^N$. Then the embedding $\widetilde{H}^K(\Omega) \hookrightarrow L^s(\Omega)$ is compact for every $s \in [2, K_*[$ if and only if $\Omega$ satisfies the assumption (D). In particular, if (D) holds, $H_0^K(\Omega)$ is compactly embedded in $L^s(\Omega)$ and the spectrum of $A$ with homogeneous Dirichlet boundary conditions in $L^2(\Omega)$ is discrete.*

*Proof.* For the first part, see [6, Theorem 2.8]. Then, since $A$ with Dirichlet boundary data is essentially selfadjoint on $C_0^\infty(\Omega)$, the discreteness of the spectrum follows by repeating the argument in [5, Theorem 4.1]. $\square$

Now, let us denote by $\Phi$ a solution (which exists by standard minimum arguments) of the following linear problem

$$\begin{cases} Au + V(x)u = 0 & \text{in } \Omega \\ \left(\frac{\partial}{\partial \nu}\right)^j u = \phi_j & \text{on } \partial\Omega \\ D^j u(x) \to 0 & \text{as } |x| \to \infty \\ j = 0, \ldots, K-1. \end{cases}$$

Notice that the regularity of $\Phi$ implies that $\Phi$ belongs to $L^t(\Omega)$ for every $t \geq 2$. Indeed, setting $\Omega_\Phi = \{x \in \Omega: \Phi(x) > 1\}$, $\Omega_\Phi$ is bounded since $\Phi$ is continuous and goes to zero at infinity and therefore $\Phi \in L^t(\Omega_\Phi)$ and

$$\int_\Omega |\Phi(x)|^t \leq \int_{\Omega_\Phi} |\Phi(x)|^t + \int_{\Omega \setminus \Omega_\Phi} |\Phi(x)|^2 < +\infty.$$

Then, the original problem $(S_{\varphi,\phi})$ can be reduced to

$$\begin{cases} Aw + V(x)w = g(x, w + \Phi) + \varphi & \text{in } \Omega \\ \left(\frac{\partial}{\partial \nu}\right)^j w = 0 & \text{on } \partial\Omega \\ D^j w(x) \to 0 & \text{as } |x| \to \infty \\ j = 0, \ldots, K-1. \end{cases}$$

More precisely, a function $u$ is a weak solution of $(S_{\varphi,\phi})$ if and only if $u \in H^K(\Omega)$, $u = w + \Phi$ and the function $w \in H_0^K(\Omega)$ satisfies

$$\forall \eta \in H_0^K(\Omega) : \int_\Omega Aw \cdot \eta + \int_\Omega V(x)w\eta = \int_\Omega g(x, w + \Phi)\eta + \int_\Omega \varphi\eta.$$

Hence, our aim is to state the existence of multiple critical points of the functional

$$I_1(u) = \frac{1}{2}\int_\Omega Au \cdot u + \frac{1}{2}\int_\Omega V(x)u^2 - \int_\Omega G(x, u + \Phi) - \int_\Omega \varphi u$$

defined on the Hilbert space $X_V = H_0^K(\Omega, V)$, endowed with the equivalent inner product

$$(u, v)_{X_V} = \int_\Omega V(x)uv + \begin{cases} \int_\Omega \Delta^m u \, \Delta^m v & \text{if } K = 2m \\ \int_\Omega \nabla\Delta^m u \, \nabla\Delta^m v & \text{if } K = 2m + 1. \end{cases}$$

We denote by $\|\cdot\|_{X_V}$ the corresponding norm induced by $(\cdot, \cdot)_{X_V}$. Following the abstract perturbation method that we will describe in the next section, let us consider the family of functionals $I_\vartheta = I(\vartheta, \cdot) : X_V \to \mathbb{R}$, $0 \le \vartheta \le 1$, defined by

$$I_\vartheta(u) = \frac{1}{2}\int_\Omega Au \cdot u + \frac{1}{2}\int_\Omega V(x)u^2 - \int_\Omega G(x, u + \vartheta\Phi) - \int_\Omega \vartheta\varphi u.$$

Standard arguments show that $I$ is a $C^1$ functional and it results

$$\frac{\partial I}{\partial \vartheta}(\vartheta, u) = -\int_\Omega g(x, u + \vartheta\Phi)\Phi - \int_\Omega \varphi u, \tag{12}$$

$$I'_\vartheta(u)[v] = \frac{\partial I}{\partial u}(\vartheta, u)[v] = \int_\Omega Au \cdot v + \int_\Omega V(x)uv - \int_\Omega g(x, u + \vartheta\Phi)v - \int_\Omega \vartheta\varphi v \tag{13}$$

for every $\vartheta \in [0, 1]$ and $u, v \in X_V$.

## 3 Proof of the results

In order to apply the method introduced by Bolle for dealing with problems with broken symmetry, let us recall the main theorem as stated in [12]. Consider two continuous functions $\varrho_1, \varrho_2 : [0, 1] \times \mathbb{R} \to \mathbb{R}$ which are Lipschitz continuous with respect to the second variable. Assume that $\varrho_1 \le \varrho_2$ and let $\psi_1, \psi_2 : [0, 1] \times \mathbb{R} \to \mathbb{R}$ be the solutions of the Cauchy problems

$$\begin{cases} \frac{\partial}{\partial \vartheta}\psi_i(\vartheta, s) = \varrho_i(\vartheta, \psi_i(\vartheta, s)) \\ \psi_i(0, s) = s. \end{cases}$$

Note that $\psi_1$ and $\psi_2$ are continuous, increasing in $s$ and $\psi_1 \leq \psi_2$. Let $J_0$ be a $C^1$–functional on a Hilbert space $X$ with norm $\|\cdot\|$. We say that a $C^1$–functional $J : [0,1] \times X \to \mathbb{R}$ is a good family of functionals starting from $J_0$ and controlled by $\varrho_1, \varrho_2$ if $J(0, \cdot) = J_0$ and if it satisfies $(\mathscr{H}_1) - (\mathscr{H}_4)$ below, where $J_\vartheta = J(\vartheta, \cdot)$.

$(\mathscr{H}_1)$ For every sequence $(\vartheta_n, u_n)$ in $[0,1] \times X$ such that

$$(J(\vartheta_n, u_n)) \text{ is bounded and } \lim_n J'_{\vartheta_n}(u_n) = 0,$$

there exists a convergent subsequence;

$(\mathscr{H}_2)$ for any $b > 0$ there exists $C_b > 0$ such that if $(\vartheta, u) \in [0, 1] \times X$ then

$$|J_\vartheta(u)| \leq b \implies \left|\frac{\partial J}{\partial \vartheta}(\vartheta, u)\right| \leq C_b(\|J'_\vartheta(u)\| + 1)(\|u\| + 1);$$

$(\mathscr{H}_3)$ for any critical point $u$ of $J_\vartheta$ we have

$$\varrho_1(\vartheta, J_\vartheta(u)) \leq \frac{\partial}{\partial \vartheta} J(\vartheta, u) \leq \varrho_2(\vartheta, J_\vartheta(u));$$

$(\mathscr{H}_4)$ for any finite dimensional subspace $W$ of $X$ it results

$$\lim_{\substack{\|u\|\to\infty \\ u \in W}} \sup_{\vartheta \in [0,1]} J(\vartheta, u) = -\infty.$$

Setting $\bar\varrho_i(s) = \sup_{\vartheta \in [0,1]} |\varrho_i(\vartheta, s)|$ we have the following result (see [12, Theorem 2.1]).

**Theorem 3.1.** *Let $\varrho_1 \leq \varrho_2$ be two velocity fields and let $\psi_1, \psi_2$ be the corresponding scalar flows. Let $J_0$ be an even $C^1$–functional on $X$ and $J$ a good family of functionals starting from $J_0$ and controlled by $\varrho_1, \varrho_2$. Let $X$ be a Hilbert decomposed as*

$$X = \bigcup_{n=0}^{\infty} X_n,$$

*where $X_0 = X_-$ is a finite dimensional subspace and $(X_n)$ is an increasing sequence of subspaces of $X$ such that $X_n = X_{n-1} \bigoplus \mathbb{R}e_n$. Consider the levels*

$$c_n = \inf_{h \in \mathcal{H}} \sup_{h(X_n)} J_0,$$

*where*

$$\mathcal{H} = \left\{ h \in C(X, X) : h \text{ is odd and } h(u) = u \text{ for } \|u\| > R \text{ for some } R > 0 \right\}.$$

*Assume that, for $n$ large, it is $c_n \geq B_1 + (B_2(n))^{\bar\beta}$ where $\bar\beta > 0$, $B_1 \in \mathbb{R}$, $B_2(n) > 0$ and*

$$\bar\varrho_i(s) \leq A_1 + A_2 |s|^{\bar\alpha}, \quad 0 \leq \bar\alpha < 1, \ A_1, A_2 \geq 0.$$

*Then $J_1$ has an unbounded sequence of critical levels if, for $n$ large, it is $(B_2(n))^{\bar\beta} > n^{\frac{1}{1-\bar\alpha}}$.*

Let us now return to our concrete framework. In order to apply the previous theorem to the functional $I(\vartheta, u)$, we need the following lemmas. In the sequel, we will denote by $c_i$ some suitable positive constants.

**Lemma 3.2.** *Let $(\vartheta_n, u_n) \subset [0,1] \times X_V$ be such that for some $C > 0$*

$$|I(\vartheta_n, u_n)| \leq C, \qquad \lim_n I'_{\vartheta_n}(u_n) = 0 \quad \text{in } X'_V.$$

*Then, up to a subsequence, $(\vartheta_n, u_n)$ converges in $[0,1] \times X_V$.*

*Proof.* Since we have $(I'_{\vartheta_n}(u_n), u_n)_{X_V} = o(\|u_n\|_{X_V})$ as $n \to +\infty$, for every $\rho \in \left]\frac{1}{\mu}, \frac{1}{2}\right[$, taking into account (2), (13) and $(g_3)$, $(g_4)$, for $n$ large it results

$$\begin{aligned} C + \rho \|u_n\|_{X_V} &\geq I_{\vartheta_n}(u_n) - \rho(I'_{\vartheta_n}(u_n), u_n)_{X_V} \\ &\geq \left(\frac{1}{2} - \rho\right) \bar{\gamma} \|u_n\|^2_{X_V} + (\rho\mu - 1) c_1 \|u_n + \vartheta\Phi\|^\mu_\mu \\ &\quad - c_2 \int_\Omega |u_n + \vartheta\,\Phi|^{p-1} |\Phi| - \int_\Omega |u_n + \vartheta\,\Phi||\varphi| - c_3 \end{aligned}$$

where $\bar{\gamma} = \min\{\gamma, 1\}$. Now, by the Young inequality, for any $\varepsilon > 0$ it results

$$\int_\Omega |u_n + \vartheta\,\Phi|^{p-1}|\Phi| \leq \varepsilon \int_\Omega |u_n + \vartheta\Phi|^\mu + \beta_{\mu,p}(\varepsilon) \int_\Omega |\Phi|^s,$$

$$\int_\Omega |u_n + \vartheta\,\Phi||\varphi| \leq \varepsilon \int_\Omega |u_n + \vartheta\,\Phi|^\mu + \beta_\mu(\varepsilon) \int_\Omega |\varphi|^{\mu'},$$

where we have set

$$\beta_\mu(\varepsilon) = \frac{\mu - 1}{\mu} \left(\frac{1}{\varepsilon\mu}\right)^{\frac{1}{\mu-1}}, \quad \beta_{\mu,p}(\varepsilon) = \frac{\mu - p + 1}{\mu} \left(\frac{p-1}{\varepsilon\mu}\right)^{\frac{p-1}{\mu-p+1}}, \quad s = \frac{\mu}{\mu - p + 1}.$$

Then, fixed $\varepsilon > 0$ sufficiently small, it follows that

$$C + \rho\|u_n\|_{X_V} \geq \left(\frac{1}{2} - \rho\right)\bar{\gamma}\|u_n\|^2_{X_V} + ((\rho\mu - 1)c_1 - (c_2 + 1)\varepsilon)\|u_n + \vartheta\Phi\|^\mu_\mu - c_4,$$

which implies the boundedness of $(u_n)$ in $X_V$. Up to a subsequence, it results $u_n \rightharpoonup u$ in $X_V$, which in view of Proposition 2.1 implies that $u_n \to u$ in $L^s(\Omega)$ for every $s \in [2, K_*[$ up to a further subsequence. Therefore, since the map

$$X_V \xrightarrow{\Upsilon} L^{\frac{K_*}{p-1}}(\Omega) \xrightarrow{(A+V(x))^{-1}} X_V, \qquad \Upsilon(u) = g(u + \vartheta\Phi)$$

is compact, a standard argument allows to prove that $u_n \to u$ in $X_V$. $\square$

**Lemma 3.3.** *For every $b > 0$ there exists $B > 0$ such that*

$$\left|\frac{\partial}{\partial \vartheta} I(\vartheta, u)\right| \leq B(1 + \|I'_\vartheta(u)\|_{X'_V})(1 + \|u\|_{X_V})$$

*for all $(\vartheta, u) \in [0,1] \times X_V$ with $|I_\vartheta(u)| \leq b$.*

*Proof.* It suffices to argue as in [13, Lemma 4.2]. $\square$

**Lemma 3.4.** *Let $\varrho_1, \varrho_2 : [0,1] \times \mathbb{R} \to \mathbb{R}$ be functions defined as*

$$-\varrho_1(\vartheta, s) = \varrho_2(\vartheta, s) = D(s^2 + 1)^{\frac{p-1}{2\mu}} \tag{14}$$

*for a suitable $D > 0$. Then*

$$\varrho_1(\vartheta, I_\vartheta(u)) \leq \frac{\partial}{\partial \vartheta} I(\vartheta, u) \leq \varrho_2(\vartheta, I_\vartheta(u))$$

*at every critical point $u$ of $I_\vartheta$. Moreover, if*

$$-\widehat{\varrho}_1(\vartheta, s) = \widehat{\varrho}_2(\vartheta, s) = D(s^2 + 1)^{\frac{1}{2\mu}} \tag{15}$$

*the same holds provided that $\phi_j = 0$ for every $j = 0, \ldots, K-1$.*

*Proof.* Arguing as in the proof of Lemma 3.2 and choosing $\rho = \frac{1}{2}$, we find $c_5 > 0$ such that

$$\|u + \vartheta \Phi\|_\mu^\mu \leq c_5 \big(I_\vartheta^2(u) + 1\big)^{\frac{1}{2}} \tag{16}$$

for every critical point $u$ of $I_\vartheta$. On the other hand by combining $(g_4)$ and (12), taking into account that $\Phi \in L^{\mu/(\mu-p+1)}(\Omega)$, we have

$$\left|\frac{\partial}{\partial \vartheta} I(\vartheta, u)\right| \leq c_6 \|u + \vartheta \Phi\|_\mu^{p-1} + c_7 \tag{17}$$

for some $c_6, c_7 > 0$ if $\phi_j \neq 0$ for some $j \in \{0, \ldots, K-1\}$ and, analogously,

$$\left|\frac{\partial}{\partial \vartheta} I(\vartheta, u)\right| \leq c_8 \|u\|_\mu \tag{18}$$

for some $c_8 > 0$ if $\phi_j = 0$ for every $j = 0, \ldots, K-1$. Therefore, putting together (16) with (17) when $\phi_j \neq 0$ and (16) with (18) when $\phi_j = 0$, the assertions follow. $\square$

Taking into account conditions $(g_2)$ and $(g_3)$, the following property can be easily shown.

**Lemma 3.5.** *For every finite dimensional subspace $W$ of $X_V$ we have*

$$\lim_{\substack{\|u\|_{X_V} \to \infty \\ u \in W}} \sup_{\vartheta \in [0,1]} I(\vartheta, u) = -\infty.$$

Let us introduce a suitable class of minimax values for the even functional $I_0$. Denote by $X_V^n$ the subspace of $X_V$ spanned by the first $n$ eigenfunctions of the operator $A + V(x)$. Let us consider

$$c_n = \inf_{h \in \mathcal{H}} \sup_{h(X_V^n)} I_0,$$

where, for a suitable constant $R > 0$, we have set

$$\mathcal{H} = \Big\{ h \in C(X_V, X_V) : h \text{ is odd and } h(u) = u \text{ for } \|u\|_{X_V} > R \Big\}.$$

Clearly, for every integer $n$, $c_n$ is a critical value of $I_0$ and $c_n \leq c_{n+1}$. Now, we need a suitable estimate on the $c'_n s$. First, let us point out that by Lemma 3.5 for all $n$ there exist $R_n > 0$ such that if $\|u\|_{X_V} > R_n$ then $I_0(u) \leq I_0(0) = 0$. Setting

$$\mathcal{H}_n = \Big\{ h \in C(D_n, X_V) : \ h \text{ is odd and } h(u) = u \text{ for } \|u\|_{X_V} = R_n \Big\},$$

where $D_n = \big\{ u \in X_V^n : \|u\|_{X_V} \leq R_n \big\}$, we deduce that

$$c_n \geq \inf_{h \in \mathcal{H}_n} \sup_{h(D_n)} I_0.$$

Arguing as in [14], we have the following result.

**Lemma 3.6.** *There exist $\bar{b} > 0$ such that, for every $n \in \mathbb{N}$,*

$$c_n \geq \bar{b} \lambda_n^{\bar{\beta}(p,K,N)}, \qquad \bar{\beta}(p,K,N) = \frac{2Kp - N(p-2)}{2K(p-2)}.$$

*Proof.* Fix $n \in \mathbb{N}$. By [14, Lemma 1.44] for every $h \in \mathcal{H}_n$ and $\rho \in \,]0, R_n[$ there exists $w$ with

$$w \in h(D_n) \cap \partial B(0,\rho) \cap (X_V^{n-1})^\perp.$$

Therefore it results

$$\max_{u \in D_n} I_0(h(u)) \geq I_0(w) \geq \inf_{u \in \partial B(0,\rho) \cap (X_V^{n-1})^\perp} I_0(u). \tag{19}$$

Notice that for every $u \in \partial B(0,\rho) \cap (X_V^{n-1})^\perp$ we have

$$I_0(u) \geq K(u), \quad \text{where } K(u) = \frac{1}{2}\|u\|_{X_V}^2 - \frac{c}{\mu}\|u\|_p^p.$$

By the Gagliardo–Nirenberg inequality there exists $c_9 > 0$ such that

$$\forall u \in X_V: \quad \|u\|_p \leq \|u\|_{K_*}^\ell \|u\|_2^{1-\ell} \leq c_9 \|u\|_{X_V}^\ell \|u\|_2^{1-\ell}, \quad \ell = \frac{N(p-2)}{2Kp}.$$

Since $u \in (X_V^{n-1})^\perp$ implies $\|u\|_2 \leq \lambda_n^{-1/2} \|u\|_{X_V}$, one obtains

$$K(u) \geq \frac{1}{2}\rho^2 - \frac{c(c_9)^p}{\mu} \lambda_n^{-(1-\ell)p/2} \rho^p.$$

Therefore, by choosing $\rho = \rho_n = c' \lambda_n^{\frac{(1-\ell)p}{2(p-2)}}$ where $c'$ is a suitable positive contant, we can assume $\rho_n < R_n$ and therefore the assertion follows by (19) and the previous inequalities. $\square$

**Remark 3.7.** As proved in [13, Lemma 5.2] for *bounded* domains (see also [21] in the case $K = 1$), the following sharp estimate holds for $n \in \mathbb{N}$

$$c_n \geq b \, n^{\frac{2pK}{(p-2)N}} \tag{20}$$

for some $b > 0$. In our framework, if $v_n$ denotes a critical point of $K(u) = \frac{1}{2}\|u\|_{X_V}^2 - \frac{c}{\mu}\|u\|_p^p$ at the level

$$b_n = \inf_{h \in \mathcal{H}_n} \sup_{h(D_n)} K(u),$$

by using suitable Morse index estimates of $v_n$, we can prove the lower estimates $b_n \geq c_{10}\|v_n\|_p^p$ and $\|v_n\|_{(p-2)N/2K} \geq c_{11} \, n^{2K/(p-2)N}$. On the other hand, if $|\Omega| = \infty$, we are not able to compare these norms and get (20). If instead $\Omega$ has finite measure, then (20) holds.

We are now ready to complete the proof of the results stated in the introduction.

**Proof of Theorem 1.1.** By combining Lemma 3.2, Lemma 3.3, (14) of Lemma 3.4, Lemma 3.5 and Lemma 3.6 the assertion follows by Theorem 3.1 with $X = X_V$, $\|\cdot\| = \|\cdot\|_{X_V}$, $X_n = X_V^n = \text{span}\{v_1, \ldots, v_n\}$ being $v_j$ the $j$–th eigenfunction of $A + V(x)$, $B(n) = \lambda_n$, $\bar{\beta} = \frac{2Kp - N(p-2)}{2K(p-2)}$ and $\bar{\alpha} = \frac{p-1}{\mu}$. If $|\Omega| < \infty$, Remark 3.7 yields the stronger conclusion.

**Proof of Theorem 1.2.** It suffices to argue as in the proof of Theorem 1.1 using (15) of Lemma 3.4 in place of (14), namely $\bar{\alpha} = \frac{1}{\mu}$

**Proof of Corollaries 1.4 and 1.5.** Taking into account Proposition 2.2, it suffices to argue as for the proof of Theorems 1.1 and 1.2 using Theorem 3.1 with $X = H_0^K(\Omega)$, $\|\cdot\| = \|\cdot\|_{K,2}$ and $X_n = \text{span}\{v_1, \ldots, v_n\}$ being $v_j$ the $j$–th eigenfunction of $A$. If $|\Omega| < \infty$, Remark 3.7 yields the stronger conclusion.

# References

[1] R.A. ADAMS, Sobolev spaces, Pure and Applied Mathematics, Vol. 65. *Academic Press* New York–London, (1975).

[2] A. BAHRI, H. BERESTYCKI, A perturbation method in critical point theory and applications, *Trans. Amer. Math. Soc.* **267** (1981), 1–32.

[3] A. BAHRI, P.L. LIONS, Morse index of some min–max critical points. I. Applications to multiplicity results, *Comm. Pure Appl. Math.* **41** (1988), 1027–1037.

[4] A. BAHRI, P.L. LIONS, Solutions of superlinear elliptic equations and their Morse indices, *Comm. Pure Appl. Math.* **45** (1992), 1205–1215.

[5] V. BENCI, D. FORTUNATO, Discreteness conditions of the spectrum of Schrödinger operators, *J. Math. Anal. Appl.* **64** (1978), 695-700.

[6] M.S. BERGER, M. SCHECHTER, Embedding theorems and quasi-linear elliptic boundary value problems for unbounded domains. *Trans. Amer. Math. Soc.* **172** (1972), 261–278.

[7] P. BOLLE, On the Bolza Problem, *J. Differential Equations* **152** (1999), 274–288.

[8] P. BOLLE, N. GHOUSSOUB, H. TEHRANI, The multiplicity of solutions in non–homogeneous boundary value problems, *Manuscripta Math.* **101** (2000), 325–350.

[9] A. CANDELA, A. SALVATORE, Multiplicity results of an elliptic equation with nonhomogeneous boundary conditions, *Topol. Methods Nonlinear Anal.* **11** (1998), 1–18.

[10] A. CANDELA, A. SALVATORE, M. SQUASSINA, Multiple solutions for semilinear elliptic systems with non–homogeneous boundary conditions, *Nonlinear Anal.* **51** (2002), 249–270.

[11] M. CLAPP, S.H. LINARES, S.E. MARTINEZ, Linking–preserving perturbations of symmetric functionals, *J. Differential Equations*, in press (2002).

[12] C. CHAMBERS, N. GHOUSSOUB, Deformation from symmetry and multiplicity of solutions in non–homogeneous problems, *Discrete Contin. Dynam. Systems* **8** (2001), 267–281.

[13] S. LANCELOTTI, A. MUSESTI, M. SQUASSINA, Infinitely many solutions for polyharmonic elliptic problems with broken symmetries, *Math. Nachr.* in press (2002).

[14] P.H. RABINOWITZ, Multiple critical points of perturbed symmetric functionals, *Trans. Amer. Math. Soc.* **272** (1982), 753–769.

[15] F. RELLICH, Über das asymptotische Verhalten der Lösungen von $\Delta u + \lambda u = 0$ in unendlichen Gebieten, *Jber. Deutsch. Math. Verein.* **53** (1943), 57–65.

[16] G.V. ROZENBLJUM, The distribution of the descrete spectrum for singular differential operators, *Soviet. Math. Dokl.* **13** (1972), 245–249.

[17] A. SALVATORE, Some multiplicity results for a superlinear elliptic problem in $\mathbb{R}^N$, *Rapp. Dip. Mat. Univ. Bari* **53** (2001).

[18] A. SALVATORE, Multiple solutions for perturbed elliptic equations in unbounded domains, *Rapp. Dip. Mat. Univ. Bari* **24** (2002).

[19] M. STRUWE, Infinitely many critical points for functionals which are not even and applications to superlinear boundary value problems, *Manuscripta Math.* **32** (1980), 335–364.

[20] M. STRUWE, Variational methods. Applications to nonlinear partial differential equations and Hamiltonian systems, Springer, Berlin, 3th edition (2000).

[21] K. TANAKA, Morse indices at critical points related to the symmetric mountain pass theorem and applications, *Comm. Partial Differential Equations* **14** (1989), 99–128.